\documentstyle[amssymb,amsfonts]{amsart}

\def\R{{\hbox{\bf R}}}
\def\T{{\hbox{\bf T}}}

\def\F{{\hbox{\rm F}}}
\def\tr{{\hbox{\rm tr}}}

\def\O{{\hbox{\bf O}}}
\def\D{{\hbox{\bf D}}}
\def\H{{\hbox{\bf H}}}

\def\Frame{{\hbox{\rm Frame}}}

\newcommand{\nabb}{{\mbox{$\nabla \mkern-13mu /$\,}}}

\def\k{{\hbox{\bf k}}}

\def\eps{\varepsilon}

\newenvironment{proof}{\noindent {\bf Proof} }{\endprf\par}
\def \endprf{\hfill  {\vrule height6pt width6pt depth0pt}\medskip}
\def\emph#1{{\it #1}}
\def\textbf#1{{\bf #1}}

\theoremstyle{plain}
  \newtheorem{theorem}[subsection]{Theorem}
  
  \newtheorem{proposition}[subsection]{Proposition}

\theoremstyle{remark}
  \newtheorem{remark}[subsection]{Remark}

\theoremstyle{definition}

\include{psfig}

\begin{document}

\title[Geometric renormalization of wave maps]{Geometric renormalization of large energy wave maps}
\author{Terence Tao}
\address{Department of Mathematics, UCLA, Los Angeles CA 90095-1555}
\email{ tao@@math.ucla.edu}
\subjclass{35J10}

\vspace{-0.3in}
\begin{abstract}
There has been much progress in recent years in understanding the existence problem for
wave maps with small critical Sobolev norm (in particular for two-dimensional wave maps with small energy); a key
aspect in that theory has been a renormalization procedure (either a geometric Coulomb gauge, or a microlocal gauge) which converts the nonlinear term into one closer to that of a semilinear wave equation.  However, both of these renormalization procedures encounter difficulty if the energy of the solution is large.  In this report we present a different renormalization, based on the harmonic map heat flow, which works for large energy wave maps from two dimensions to hyperbolic spaces.
We also observe an intriguing estimate of ``non-concentration'' type, which asserts roughly speaking that if the energy of a wave map concentrates at a point, then it becomes asymptotically self-similar.
\end{abstract}

\maketitle

\section{Introduction}

Let $n,m \geq 1$ be integers, and let $\R^{1+n}$ be Minkowski space $\{ (t,x): t \in \R, x \in \R^n \}$
with the usual metric $g_{\alpha \beta} x^\alpha x^\beta = -dt^2 + dx^2$.  We parameterize this space by Greek indices
$\alpha, \beta = 0,1,\ldots,n$, raised and lowered in the usual manner.

Let $(N,h)$ be a smooth connected complete $m$-dimensional Riemannian manifold without boundary, and let $I \subseteq \R$ be an interval. A \emph{Schwartz wave map}  from a slab $I \times \R^n$ in Minkowski space to $N$ is defined to be any smooth map 
$\phi: I \times \R^n \to N$ which decays rapidly to some constant $\phi(\infty) \in N$ at spatial infinity, and whose derivatives all decay rapidly to 0.\footnote{In other words,
for each fixed $t$, $\phi(t,x)$ converges to $\phi(\infty)$ faster than any negative power of $|x|$, and all derivatives of $\phi$ converge to zero faster than any negative power of $|x|$.  We need to consider Schwartz wave maps rather
than classical wave maps (which are constant outside of in a finite union of light cones) because we will need to take
(nonlinear) Littlewood-Paley projections of these wave maps, obtained via the heat flow for harmonic maps,
in our argument.}, and is a (formal) 
critical point of the 
Lagrangian
\begin{equation}\label{lagrangian}
 \int_{\R^{1+n}} \langle \partial^\alpha \phi(t,x), \partial_\alpha \phi(t,x) \rangle_{h(\phi(t,x))}\ dt dx.
\end{equation}
The Euler-Lagrange equation for Schwartz wave maps can then be written as
\begin{equation}\label{cov}
 (\phi^*\nabla)^\alpha \partial_\alpha \phi = 0,
\end{equation}
where $(\phi^*\nabla)^\alpha$ is covariant differentiation on the vector bundle $\phi^*(TN)$ 
with respect to the pull-back $\phi^*\nabla$ via $\phi$ of the Levi-Civita connection $\nabla$ on $N$.
Wave maps are of course the Minkowski space analogue of harmonic maps.

In this paper we shall restrict primarily to the case when $n=2$, $m \geq 2$ (the $m=1$ case being reducible to
the free (linear) wave equation) and when $N$ has constant negative curvature $\kappa < 0$, thus
\begin{equation}\label{constant-curvature}
 [\nabla_X,\nabla_Y] Z - \nabla_{[X,Y]} Z =: R(X,Y) Z = \kappa ( \langle Y,Z \rangle_h X - \langle X,Z \rangle_h Y )
\end{equation}
for all smooth vector fields $X,Y,Z$ on $N$.  Note that $\kappa$ can be normalized to $-1$ by scaling the metric
$h$.  A typical example of such a manifold $N$ (with $\kappa = -1$) 
is the hyperbolic space $\H^m := O(m,1)/O(m) = SO(m,1)/SO(m)$, where $O(m)$ is the orthogonal group on $\R^m$ 
and $O(m,1)$ is the Lorentz group on $\R^{1+m}$, with $O(m,1)$ and $SO(m,1)$ being the respective orientation-preserving components of these groups.  More generally, one can quotient $\H^m$ by
a discrete subgroup of $O(m,1)$ (acting on the left in the usual manner), and as is well known one can recover all constant negative curvature $m$-dimensional manifolds $N$ this way, up to isometry and normalization of $\kappa$.  The case $N = \H^2$ is of particular interest as it arises as a model problem for the Einstein equations under some additional symmetry assumptions (see e.g. \cite{bgm} for a discussion).

\begin{remark} With the constant negative curvature assumption, Schwartz maps $\phi$ become topologically trivial (contractable to a constant map); by compactifying infinity, this remark is equivalent to the assertion that any map from a sphere $S^n$ to $N$ is topologically trivial.  When $N$ is hyperbolic space $\H^m$ this claim is clear since $N$ is contractible; more generally, if $N$ is $\H^m$ quotiented out by a discrete group, then the claim follows by lifting the map up to hyperbolic space (the discreteness of the fibers keeps the monodromy trivial).  Later on we will see that these maps can in fact be explicitly deformed into a constant map by means of the harmonic map heat flow.  Note that while $N$ itself need not be parallelizable or even orientable, this will not particularly concern us since maps into $N$ can be lifted up into $\H^m$, which is both parallelizable and orientable.  Indeed one could reduce to the case $N=\H^m$ by
lifting the wave map up to hyperbolic space, and observing that the wave map equation \eqref{cov} is preserved by this lift.
\end{remark}

In local co-ordinates $\phi = \phi^a$ on $N$, the wave map equation \eqref{cov} can be written as 
\begin{equation}\label{wm-coord}
 \partial^\alpha \partial_\alpha \phi^a = - \Gamma^a_{bc}(\phi) \partial_\alpha \phi^b \partial^\alpha \phi^c
\end{equation}
where the $\Gamma^a_{bc}$ are Christoffel symbols; in particular we see that $\phi$ evolves by a non-linear 
wave equation.  It is then natural
to consider the Cauchy problem, specifying initial data $\phi[t_0] := (\phi(t_0), \partial_t \phi(t_0))$ 
at some initial time $t_0 \in \R$.  Of course, this initial data $\phi[t_0]$
must take values in the tangent bundle $TN$ of $N$,
and decay rapidly to $(\phi(\infty),0)$ as $x \to \infty$; we refer to such data as \emph{Schwartz initial data}.  

The Cauchy problem for wave maps has been extensively studied (see references); we refer the 
interested reader to the surveys in \cite{kman.barrett}, \cite{kman.selberg:survey}, \cite{shatah-struwe}, \cite{struwe.barrett}, \cite{tataru:survey}.  For arbitrary dimension $n$ and arbitrary targets $N$, one can easily show by general energy methods (which indeed apply to any nonlinear wave equation, see e.g. \cite{sogge:wave}) that for any Schwartz\footnote{The smoothness condition has been relaxed substantially, see for instance \cite{kl-mac:null3}, \cite{kman.selberg}, \cite{keel:wavemap}, \cite{tataru:wave1}, \cite{tataru:wave2}, \cite{tataru:wave3}, \cite{tataru:survey} and the discussion below; however we shall restrict our attention here to the Schwartz category.  The rapid decay assumptions can also be removed by finite speed of propagation, though possibly at the cost of creating a domain of existence which is not a spacetime slab $I \times \R^n$.} initial data $\phi[t_0]$, one has a unique \emph{local-in-time} Schwartz solution to \eqref{cov} in some spacetime slab $I \times \R^n$, where $I$ is a neighborhood of $t_0$, attaining the value of $\phi(\infty)$ at infinity (thus the value of $\phi$ at spatial infinity is fixed by the wave map evolution).  Indeed, energy methods can be used to show such smooth solutions can be continued
whenever the norm $\| \nabla_{t,x} \phi(t) \|_{H^{n/2+\eps}(\R^n)}$ remains bounded for some $\eps > 0$.
In particular, by use of the conserved energy
\begin{equation}\label{energy-cov}
 E(\phi) = E(\phi[t]) := \int_{\R^n} \frac{1}{2} |\nabla_{t,x} \phi(t,x)|_h^2\ dx
\end{equation}
we see that one has global existence of Schwartz solutions in one dimension; see
\cite{gu}, \cite{lady}, \cite{ginebre}, \cite{shatah}, \cite{keel:wavemap}.  The one dimensional wave map equation
is actually rather special, being completely integrable and thus supporting special solutions such as homoclinic breather
solutions; see \cite{strauss}.

Now we turn to higher dimensions $n \geq 2$.  Since the wave map equation \eqref{cov}
is invariant under the scaling
\begin{equation}\label{cov-scaling}
\phi(t,x) \leftarrow \phi(\frac{t}{\lambda}, \frac{x}{\lambda})
\end{equation}
for any $\lambda > 0$, we see that the two-dimensional case $n=2$ is critical with respect to the energy \eqref{energy-cov}, as this energy is now scale-invariant (dimensionless) in this case.  For higher dimensions $n > 2$ it is known
that singularities can form\footnote{Despite the presence of these singularities, it is known that global 
\emph{weak} solutions exist globally in time, see \cite{muller}, \cite{fms:compact}, although it is quite possible that uniqueness and energy conservation will fail for these solutions.}
 in finite time when the manifold $N$ is positively curved, 
\cite{shatah}, \cite{shatah.shadi.blow}, \cite{css}, \cite{shatah-struwe} (and there is strong numerical evidence to suggest that there is also blowup when $n=2$ \cite{isenberg}, \cite{bct} in this case); when $n \geq 7$ one can also 
obtain singularities for negatively curved manifolds \cite{css}, despite such manifolds being well-behaved for other equations (such as the harmonic map heat flow).  However, one still has global solutions for small data.  This was 
first observed for initial data sufficiently close to a point \cite{chq} or a geodesic \cite{sideris.harmonic}, and for data which is small in the critical Besov norm $\| \nabla_{t,x} \phi \|_{\dot B^{n/2-1}_1(\R^n)}$ \cite{tataru:wave1}, \cite{tataru:wave2},
assuming that the manifold $N$ is uniformly smooth.  This was achieved by viewing the wave map equation as
a nonlinear wave equation in $\phi$ (as in \eqref{wm-coord}) and iterating in a certain sophisticated
function space.  These Besov space results appear to be essentially the limit of the method; for instance, it is known that such 
iteration methods cannot work if the data is merely assumed to be small in $\| \nabla_{t,x} \phi \|_{\dot H^{n/2-1}(\R^n)}$ (see \cite{keel:wavemap}, \cite{nak:wavemap}, \cite{dag}).  However, this small Sobolev norm global regularity result turned out to be achievable upon either a microlocal renormalization of the wave map equation, or by passing to a derivative formulation and selecting a special gauge (such as the Coulomb gauge).  In the case of the sphere 
$N = S^{m-1}$, this was achieved for $n \geq 2$ in \cite{tao:wavemap}, \cite{tao:wavemap2} via the microlocal renormalization method; for the more general class of boundedly parallelizable manifolds (which includes the hyperbolic spaces $\H^m$), this was achieved for $n \geq 5$ in \cite{kr:wavemap} (by a hybrid of the microlocal renormalization and Coulomb gauge
methods) and then for $n \geq 4$ in
\cite{shsw:wavemap}, \cite{nahmod} and $n \geq 3$ in \cite{krieger:3d} (by use of the Coulomb gauge).  For the hyperbolic
plane $N = \H^2$ the $n=2$ case was recently treated in \cite{krieger:2d} (again using the Coulomb gauge), while
for manifolds which are uniformly isometrically embeddable in Euclidean space the result was obtained for all $n \geq 2$
in \cite{tataru:wave3} (together with more precise well-posedness results).  While the result in \cite{tataru:wave3} does not directly cover the hyperbolic spaces $\H^m$, which cannot be uniformly isometrically embedded, it may be possible that the argument can be modified to treat this case by first quotienting $\H^m$ by a discrete group in order to compactify the target, and lifting back to $\H^m$ at the end of the argument.

Thus we now have a fairly satisfactory regularity theory when the critical Sobolev norm is small; in particular in 
two dimensions we have a small energy regularity theory for a reasonably large class of target manifolds.  Now we turn to the question of large energy data 
at $n=2$.  In light of the  numerical work in \cite{isenberg}, \cite{bct} we do not expect global regularity in 
the positive curvature case, but it was still conjectured (see e.g. \cite{klainerman:unified}) that one has 
global regularity for large energy data in negative curvature manifolds such as $\H^m$; note that this would in fact imply an unconditional global regularity result for smooth data thanks to finite speed of propagation.  In the case where the wave map has some $U(1)$ symmetry (either rotation-invariant, or rotation equivariant), such results have already been obtained in
\cite{christ.spherical.wave} (for spherically symmetric maps) and \cite{shatah.shadi.eqvt}, 
\cite{shsw:wavemap-eqvt}, \cite{shatah-struwe} (for equivariant maps).  The spherically symmetric regularity result 
in fact extends to arbitrary targets \cite{struwe:radial-sphere}, \cite{struwe:radial-general} but the equivariant one seems restricted to sufficiently
``negatively curved'' targets; see \cite{struwe:equivariant} (and the numerics in \cite{bct}, \cite{isenberg}) for 
further discussion of this issue, which is in particular related to the existence of non-constant \emph{harmonic} 
maps from $S^2$ to $N$, which among other things can be used to generate non-trivial finite-energy stationary solutions 
to \eqref{cov}).  
numerics that strongly suggest blowup occurs for instance when the target is a sphere, even for equivariant data.
In \cite{struwe:equivariant} it is shown that for equivariant data, blowup can only occur if the wave map relaxes
to a (rescaled) harmonic map from $\R^2$ to $S^2$ for a sequence of times converging to the blowup time.

The conjecture of large energy global regularity for wave maps into $\H^m$, if true, would be
analogous to other critical large energy global regularity results in the literature, notably that of the 
energy-critical nonlinear wave equation
$\Box u = u^5$ in three dimensions (see e.g. \cite{shatah-struwe}), the critical Yang-Mills-Higgs equation 
in three dimensions \cite{keel},  and the energy-critical nonlinear Schr\"odinger
equation $iu_t + \Delta u = |u|^4 u$ in three dimensions (see \cite{grillakis:scatter}, \cite{borg:book} for the radial case,
and \cite{gopher} for the general case).  These large energy results are substantially more difficult than
the small energy theory, and require at least two new ingredients:

\begin{itemize}

\item A perturbation, local existence, and theory suitable for large energy solutions, assuming that a suitable spacetime norm (e.g.
the $L^{10}_{t,x}$ norm for the energy-critical three-dimensional nonlinear Schr\"odinger equation, or $L^4_t L^{12}_x$ norm for the energy-critical three-dimensional nonlinear wave equation or Yang-Mills-Higgs equation) is bounded;

\item Some sort of energy non-concentration argument which establishes that such a spacetime norm is bounded.

\end{itemize}

The energy non-concentration argument generally relies on some sort of monotonicity formula, often arising
from a Morawetz inequality, although in some cases (notably that of the energy-critical Schr\"odinger equation)
such an inequality is not sufficient by itself to control the desired spacetime quantity; one must supplement that
inequality with additional arguments, such as the induction-on-energy strategy pioneered by Bourgain \cite{borg:book}
and which also plays a fundamental role in \cite{gopher}.  On the other hand, in these types of arguments,
the large energy perturbation analysis does not require such a Morawetz inequality, and instead proceeds by
a variant of the standard existence theory, using Strichartz estimates as the primary tool.  In order to do this,
one relies crucially on the \emph{fungibility in time} of spacetime norms such as the $L^{10}_{t,x}$; what we mean 
by this is that if the $L^{10}_{t,x}$ norm is large but bounded, then one can decompose the time interval into 
a bounded number of sub-intervals, such that the $L^{10}_{t,x}$ norm is small on each interval.  This allows 
one (at least in principle) to obtain a perturbation theory for large $L^{10}_{t,x}$ norm solutions from that for 
small $L^{10}_{t,x}$ norm solutions.  This is in contrast with norms such as $L^\infty_t L^6_x$ which are clearly
bounded by the energy for (say) the energy-critical NLS, but are not fungible in time and so cannot be used directly
to establish a global regularity result for this equation.

We do not fully resolve either of these two issues here, however we can announce some progress on both.  More
precisely, we are able to establish a large energy perturbation theory for wave maps, in the case when the target is a space of constant negative curvature, although the actual details of this theory are rather technical and will be only partially sketched here.  Secondly, we have a simple monotonicity formula argument which shows that
if the energy of a wave map (with arbitrary target) concentrates at a point, then it must be asymptotically self-similar
in an averaged sense (see Proposition \ref{alss}).

We now discuss the gauge renormalization (leaving the discussion of the non-concentration result to Section
\ref{alss-sec}).  For the model problem of constant negative curvature, a direct attempt to mimic either the microlocal renormalization procedure in
\cite{tao:wavemap}, \cite{tao:wavemap2}, \cite{kr:wavemap}, \cite{tataru:wave3} or the Coulomb gauge
renormalization procedure in \cite{shsw:wavemap}, \cite{nahmod}, \cite{krieger:3d}, \cite{krieger:2d} runs into
difficulty (although the Coulomb gauge renormalization will work in the abelian case, in which the target is
a two-dimensional space such as $\H^2$, so that the gauge group $SO(2)$ is abelian).  More specifically, one begins to encounter difficulty in large energy
in keeping the microlocal gauge change approximately unitary, and in the Coulomb gauge one has problems establishing
uniqueness, regularity, and ellipticity of the gauge in the non-abelian case where the target has dimension $m>2$ and
so the gauge group $SO(m)$ is non-abelian.  Also, the microlocal gauge has not yet been successfully applied to
manifolds such as hyperbolic space, which are not easy to embed inside Euclidean spaces; meanwhile,
the Coulomb gauge introduces some additional (and somewhat
artificial) singularities at the spatial frequency origin $\xi = 0$ (arising from the Hodge decomposition)
which cause additional technical complications, especially in two dimensions (see \cite{krieger:3d}, \cite{krieger:2d}
for a demonstration of this phenomenon).  Finally, while both these gauges are well suited for establishing regularity of
a \emph{single} solution, the task of establishing a perturbation theory - i.e. comparing two nearby solutions - is also
not trivial, as one must then establish that the gauge transforms themselves have suitable continuity properties, which is an extremely delicate issue in the critical setting.

In this report we give a new renormalization procedure, which is based on the harmonic map heat flow, which is
intrinsic to the geometry and seems a extremely natural choice of gauge; we shall dub it the ``canonical heat-temporal gauge'' or ``caloric gauge''.  This caloric gauge has essentially the same effect
as the microlocal gauge change (which can be viewed as a discretized, partially linearized, version of
the harmonic map heat flow renormalization), in that it replaces the nonlinearity by a ``paraproduct'' variant
in which at least one derivative falls on a low frequency term.  However, unlike the microlocal gauge change, the caloric gauge works well in large energies and negative curvature targets, essentially thanks to the famous result of Eells and
Sampson \cite{eells} which shows that the harmonic map heat flow for those targets is globally smooth and converges
to a constant.  The caloric gauge also is well adapted to perturbation theory, since the harmonic map heat flow also enjoys a good perturbation theory.  When compared against the Coulomb gauge, the caloric gauge has the advantage of
renormalizing the nonlinearity into a slightly better form (with no singularity near the spatial frequency origin);
we shall give a heuristic comparison of the two renormalized nonlinearities using Littlewood-Paley analysis later
on in this report.  It is also insensitive to non-abelian behaviour in the gauge group, and does not develop difficulties
with uniqueness or regularity, as long as the target is negatively curved.

The discussion here will consist of rigorous geometric computations, and non-rigorous analytical heuristics; a rigorous
treatment of the perturbation theory which is informally discussed here will appear elsewhere (and in a much lengthier paper).

\emph{Acknowledgements.} The author is indebted to Ben Andrews and Andrew Hassell for a crash course in Riemannian geometry and manifold embedding, and in particular to Ben Andrews for explaining the harmonic map heat flow.  The author also thanks Mark Keel for background material on wave maps, and Daniel Tataru for sharing some valuable insights on multilinear estimates and function spaces.  The author thanks Andrew Hassell and the Australian National University for their hospitality when much of this work was conducted.  The author is a Clay Prize Fellow and is supported by a grant from the Packard Foundation.

\section{The derivative formulation}

Here and in the rest of the paper we fix $N$ to be a smooth complete Riemannian manifold of constant negative curvature $\kappa < 0$.
We now describe the (standard) derivative formulation of the wave map equation \eqref{cov}, writing $\nabla_{t,x} \phi$ in terms
of an orthonormal frame $e(t,x)$; this formulation was used for instance in
\cite{kr:wavemap}, \cite{shsw:wavemap}, \cite{nahmod} (or in \cite{helein} for harmonic maps).  Our discussions here shall
be primarily geometric rather than analytic.  As we shall need
this derivative formulation not just for the wave map equation, but also for the harmonic map heat flow (and in fact for
a concatenation of the two equations), we shall first describe the derivative formulation of a general smooth 
map $\phi: \Omega \to N$, where $\Omega$ is an open contractible subset of either a Minkowski space $\R^{1+n}$ or a Euclidean space $\R^n$.
In either case we will use the co-ordinates $x^\alpha$ to parameterize the domain $\Omega$, with the associated
partial derivative operators $\partial_\alpha$.  For this general discussion we will not need the metric structure of $\Omega$
(i.e. we will not raise and lower indices), although we will use the fact that $\Omega$ is flat in the sense that the 
partial derivatives $\partial_\alpha$ commute with each other, thus $[\partial_\alpha, \partial_\beta] = 0$.

We recall the \emph{orthonormal frame bundle} $\Frame(N)$ of $N$, defined as the collection of all pairs
$(\phi, e)$ where $\phi \in N$ and $e: \R^m \to T_\phi N$ is an orthogonal transformation (or equivalently, an orthonormal basis of $T_\phi N$).  It is geometrically obvious that, $\Frame(N)$ is a smooth principal $O(m)$-bundle of $N$, where $O(m)$ is the orthogonal group on $\R^m$ (not necessarily orientation preserving).  Since $\Omega$ is contractible, we can lift any map $\phi: \Omega \to N$ to a map $(\phi,e): \Omega \to \Frame(N)$ into 
the orthonormal frame bundle.  In fact there are multiple such lifts; given any
smooth function $U: \Omega \to O(m)$, we have the \emph{gauge transform} $(\phi,e) \mapsto (\phi, eU^{-1})$ which
transforms any lift of $\phi$ to any other lift of $\phi$; indeed, since $\Frame(N)$ is a principal $O(m)$-bundle,
all such lifts can be obtained from a single specified lift by a (unique) gauge transform.

Let $(\phi,e): \Omega \to \Frame(N)$ be a smooth map.  The bundle $\phi^*(TN)$ (the pullback of the tangent bundle of $N$
under $\phi$) is an $m$-dimensional vector bundle on $\Omega$, and the inverse $e^{-1}$ of the
orthonormal frame $e$ can be used to identify this with the trivial bundle $\Omega \times \R^m$.
In particular, since the derivatives $\partial_\alpha \phi$ are sections of $\phi^*(TN)$, we can define the $\R^m$-valued one-form $\psi_\alpha: U \to \R^m$ by this procedure, thus
$$ \psi_\alpha := e^{-1} \partial_\alpha \phi$$
or equivalently
\begin{equation}\label{psi-def}
 \partial_\alpha \phi = e \psi_\alpha
\end{equation}

The Levi-Civita connection $\nabla$ on $TN$ induces a pullback connection $\phi^*(\nabla)$ on $\phi^*(TN)$, which in turn
induces a connection on the bundle $\Omega \times \R^m$ by the above identification.  This connection is compatible with the gauge group $O(m)$ (since the original Levi-Civita connection is also), and is given by the covariant derivatives
\begin{equation}\label{cov-deriv-def}
\D_\alpha := \partial_\alpha + A_\alpha.
\end{equation}
where the $A_\alpha: \Omega \to o(m)$ take values in $o(m)$, the skew-symmetric linear transformations of $\R^m$ (and
the Lie algebra of $O(m)$) and are defined by
\begin{equation}\label{A-def}
 (\phi^*\nabla)_\alpha e = e A_\alpha.
\end{equation}
From the Leibnitz rule and \eqref{A-def} we have
\begin{equation}\label{cov-cov}
 (\phi^*\nabla)_\alpha (ef) = e(\D_\alpha f)
\end{equation}
for any smooth scalar function $f: \Omega \to \R^m$.  Since the Levi-Civita connection $\nabla$ is torsion free, we have
$$ \nabla_{\partial_\alpha \phi} \partial_\beta \phi = \nabla_{\partial_\beta \phi} \partial_\alpha \phi$$
which upon pulling back by $\phi$ becomes
$$ (\phi^* \nabla)_\alpha \partial_\beta \phi = (\phi^* \nabla)_\beta \phi \partial_\alpha \phi.$$
From \eqref{psi-def}, \eqref{cov-cov} we thus have the fundamental \emph{zero torsion identity}
\begin{equation}\label{zero-torsion}
\D_\alpha \psi_\beta - \D_\beta \psi_\alpha = 0.
\end{equation}
This identity allows us (at least in principle) to write any of the fields $\psi_\alpha$ in terms of a single 
field $\psi_\beta$, thus converting the derivative wave map system into something resembling a scalar equation.

\begin{remark} The tangent bundle $T \Frame(N)$ of $\Frame(N)$ can be canonically identified with $\Frame(N) \times \R^m \times o(m)$, since for any $(\phi,e) \in \Frame(N)$, the orthonormal basis $e$ can be used to parameterize the horizontal space $T_\phi N$ (lifted up to $T_{(\phi,e)} \Frame(N)$ via the Levi-Civita connection) as $\R^m$, while the vertical space of the bundle $\Frame(N)$ at $(\phi,e)$ can be parameterized as $o(m)$.  The fields $(\psi_\alpha,A_\alpha)$ can thus be thought of as the derivatives $\partial_\alpha (\phi,e)$ of the original fields $(\phi,e)$ with respect to this identification.  Furthermore, one can show that the Lie algebra structure of
this $\R^m \times o(m)$-bundle induces an action of the special Lorentz group $O(m,1)$, which then acts transitively on $\Frame(N)$.  This in turn can be used to verify the assertion made earlier that $N$ is the quotient of $\H^m$ by a discrete group.  We were however unable to exploit this group structure to any particular advantage, and shall mostly treat the fields $\phi$, $e$ and $\psi_\alpha,A_\alpha$ separately.
\end{remark}

The connection $A_\alpha$ has a \emph{curvature tensor} $F_{\alpha \beta}: \Omega \to o(m)$ defined by
\begin{equation}\label{curv-def}
F_{\alpha \beta} := [\D_\alpha,\D_\beta] = \partial_\alpha A_\beta - \partial_\beta A_\alpha + [A_\alpha, A_\beta].
\end{equation}
We can use the constant negative curvature hypothesis \eqref{constant-curvature} to compute $F_{\alpha \beta}$. Indeed, from \eqref{constant-curvature} we have
$$ [(\phi^*\nabla)_\alpha, (\phi^*\nabla)_\beta] Z = 
\kappa (\langle \partial_\beta \phi, Z \rangle_h \partial_\alpha \phi - \langle \partial_\alpha \phi, Z \rangle_h \partial_\beta \phi)$$
for all sections $Z$ of $\phi^*(TN)$; applying \eqref{psi-def}, \eqref{cov-cov}, \eqref{curv-def} we obtain the \emph{curvature identity}
\begin{equation}\label{curvature-identity}
F_{\alpha \beta} = \kappa \psi_\alpha \wedge \psi_\beta
\end{equation}
where we use $u \wedge v \in o(m)$ to denote the skew-symmetric linear operator on $\R^m$ defined by\footnote{The cancellation inherent in this wedge product operation will be largely unexploited, in contrast to \cite{krieger:2d}, \cite{krieger:3d}, where
this type of ``$Q_{ij}$ null structure'' is needed to compensate for the singularity of $\nabla^{-1}$ in the Coulomb 
gauge.  We will, however, rely crucially on the negativity of $\kappa$ to ensure that the heat flow converges properly.}
$$ (u \wedge v) w := u \langle v, w \rangle - v \langle u, w \rangle.$$
This identity will allow us (given a suitable fixing of gauge, see below) 
to recover $A$ as some sort of quadratic integral of $\psi$.

\begin{remark} In the Riemann surface case $m=2$, the gauge group $O(m)$ is abelian and so the commutator term $[A_\alpha,A_\beta]$ vanishes in \eqref{curv-def}.  This is in particular useful for placing $A$ in the Coulomb gauge.
However in the canonical heat-temporal gauge, we will not be unduly bothered by this commutator term, indeed in that gauge this commutator will have a number of vanishing components and will thus essentially be ignored.
\end{remark}

We refer to $(\phi,e)$ as the \emph{undifferentiated fields} and $(\psi_\alpha,A_\alpha)$ as the \emph{differentiated fields}.  The gauge freedom $e \mapsto eU^{-1}$ mentioned earlier affects these fields by the formulae
\begin{equation}\label{gauge-def}
(\phi,e) \mapsto (\phi, eU^{-1}); \quad 
(\psi_\alpha, A_\alpha) \mapsto (U \psi_\alpha, U A_\alpha U^{-1} + (\partial_\alpha U) U^{-1}); 
\quad F_{\alpha \beta} \to U F_{\alpha \beta} U^{-1}.
\end{equation}
In particular we observe that the magnitudes $|\psi_\alpha|$ and $|F_{\alpha \beta}|$ are independent of the
choice of gauge.  Here we use $|F_{\alpha \beta}|$ to denote the Hilbert-Schmidt norm, thus
$|F| := \tr( FF^* )^{1/2}$.

\begin{remark}\label{converse}
Clearly from \eqref{psi-def} and \eqref{A-def} we can obtain the differentiated fields from the undifferentiated fields.  Now we consider the reverse question, whether one can obtain the undifferentiated
fields from the differentiated fields.  To fix boundary conditions, let us now assume that $\Omega$ is  
convex, unbounded, and connected at infinity, and that $(\phi,e)$ converges 
rapidly (in the Schwartz sense, i.e. all derivatives $\nabla_{x,t}^N (\phi,e)$
decay faster than any polynomial) to some point $(\phi(\infty),e(\infty)) \in \Frame(N)$ at infinity\footnote{We remark that if $\phi$ decays rapidly to $\phi(\infty)$ at infinity and $e(\infty)$ is any orthonormal frame of $T_{\phi(\infty)} N$, then one can easily find a lift $(\phi,e): \Omega \to \Frame(N)$ of $\phi: \Omega \to N$ which converges to $(\phi(\infty),e(\infty))$ at infinity.  If $N$ is hyperbolic space $\H^m$ this is clear since $\H^m$ is parallelizable and one can simply pull back an orthonormal frame on $\H^m$ which equals $e(\infty)$ at $\phi(\infty)$; more generally one can lift the map $\phi: \Omega \to N$ to a map $\tilde \phi: \Omega \to \H^m$ (retaining the Schwartz property at infinity, obtain a frame $(\tilde \phi,\tilde e)$ for the lift, and then push the frame back down to $N$.} This then implies that $(\psi_\alpha, A_\alpha)$ are Schwartz. To reconstruct $(\phi,e)$ from $(\psi_\alpha,A_\alpha)$ one simply solves the ODE \eqref{psi-def}, 
\eqref{A-def} from the given data at infinity. A simple calculation using \eqref{zero-torsion}, \eqref{curvature-identity}, and Frobenius's theorem\footnote{Alternatively, if $0 \in \Omega$, one can use a radial gauge $x^\alpha (\phi^* \nabla)_\alpha e = 0$) and propogate \eqref{psi-def}, \eqref{A-def} along radial directions; one can then use \eqref{zero-torsion}, \eqref{curvature-identity} and a Gronwall inequality argument to show that $(\psi_\alpha,A_\alpha)$ are indeed the derivative fields of $(\phi,e)$.  Again, we omit the details.} shows that this system is indeed solvable; we omit 
the details.  
\end{remark}

To summarize, we now have have two equivalent ways to view a map $(\phi,e):\Omega \to \Frame(N)$;
either in an undifferentiated formulation, in which $(\phi,e)$ can move freely (but do not take values in
a linear space), or in a differentiated formulation, in which the fields $(\psi_\alpha, A_\alpha)$ now live 
in linear spaces but are constrained to obey the compatiblity conditions 
\eqref{zero-torsion}, \eqref{curvature-identity}.  In the Schwartz category
these two viewpoints are equivalent as long as we fix $(\phi(\infty), e(\infty))$ in advance.  As we shall be needing
both the differentiated and undifferentiated fields at various stages our analysis, we shall however not take advantage of this redundancy, and instead consider all the fields simultaneously.

\subsection{The derivative formulation of the wave map equation}

We now specialize to Schwartz wave maps $\phi: I \times \R^n \to N$ for some compact interval $I$, which equal some
fixed point $\phi(\infty) \in N$ at infinity.  Given any fixed frame $e(\infty)$ in $T_{\phi(\infty)} N$, we can
lift $\phi$ to a map $(\phi,e): I \times \R^n \to \Frame(N)$ which converge rapidly to $(\phi(\infty),e(\infty))$ at infinity.

The wave map equation \eqref{cov} can be rewritten using \eqref{psi-def}, \eqref{cov-cov} as
$$ \D^\alpha \psi_\alpha = 0$$
where we of course now use the Minkowski metric to raise and lower indices.  Combining this equation with the
compatibility conditions \eqref{zero-torsion}, \eqref{curvature-identity} and the definitions \eqref{psi-def}, \eqref{A-def} we can thus write the wave map equation as the system (cf. \cite{shsw:wavemap}, \cite{nahmod}, \cite{krieger:3d}, \cite{krieger:2d})
\begin{equation}\label{wmeq}
\begin{split}
\partial_\alpha \phi &= e \psi_\alpha\\
(\phi^*\nabla)_\alpha e &= e A_\alpha\\
\D^\alpha \psi_\alpha &= 0\\
\D_\alpha \psi_\beta - \D_\beta \psi_\alpha &= 0 \\
F_{\alpha \beta} &= \kappa \psi_\alpha \wedge \psi_\beta.
\end{split}
\end{equation}
Note that these equations are both over-determined (the compatibility conditions \eqref{zero-torsion}, \eqref{curvature-identity} enforcing a number of constraints
between the fields) and under-determined (since we still have the gauge freedom \eqref{gauge-def}).  In light of
Remark \ref{converse} it is possible to omit the first two equations, and just consider the last three equations
as a system for the differentiated fields $\psi_\alpha$, $A_\alpha$; since can then recover the undifferentiated fields
$\phi$, $e$ using Frobenius's theorem (and the boundary data $\phi(\infty), e(\infty)$) as discussed earlier; note that we have relied crucially on the constant curvature in order to forget about $\phi$ and $e$.  However, 
we will retain both the differentiated and
undifferentiated fields in our system, as the presence of the undifferentiated fields makes the compatibility conditions
automatic.

We refer to the quantity $u := \D^\alpha \psi_\alpha$ as the \emph{wave-tension field} (as opposed to the \emph{heat-tension} field $\psi_s = \D_k \psi_k$ which will appear later); thus wave maps are precisely those maps whose
wave-tension field vanishes.

The system \eqref{wmeq} is not yet well-posed in time, even locally, because of the gauge freedom \eqref{gauge-def}.
There are a number of options available to fix the gauge freedom and recover well-posedness; examples of such gauges include the \emph{temporal gauge} $A_0 = 0$, the \emph{Lorenz gauge} $\partial_\alpha A^\alpha = 0$, and the \emph{Coulomb gauge} $\partial_j A_j = 0$, where Latin indices such as $j$ will run
 over the spatial indices $1, \ldots, n$ with the usual summation conventions.  The latter gauge has been exploited
for instance in \cite{nahmod}, \cite{shsw:wavemap}, \cite{krieger:3d}, \cite{krieger:2d}; one advantage of
this gauge is that it makes $A$ a quadratic function of $\psi$.
Indeed, from Hodge theory one can easily verify in the Coulomb gauge that
\begin{equation}\label{hodge}
\begin{split}
A_\alpha &= \Delta^{-1} (\partial_j \partial_j A_\alpha)\\
&= \Delta^{-1}( \partial_j (F_{j\alpha} + [A_\alpha,A_j]) + \partial_\alpha \partial_j A_j)\\
&= \Delta^{-1} \partial_j (\kappa \psi_j \wedge \psi_\alpha + [A_\alpha,A_j]),
\end{split}
\end{equation}
where $\Delta$ is the spatial Laplacian.  Schematically, we thus have $A = \nabla_x^{-1}(\psi \psi + A A)$.  This 
identity allows the wave map equation to be quite tractable in high dimensions 
(see e.g. \cite{nahmod}, \cite{shsw:wavemap} for the $n \geq 4$ cases), however the inverse derivative $\nabla_x^{-1}$ causes some difficulties in low dimensions, although they can eventually be surmounted (see \cite{krieger:3d}, \cite{krieger:2d}).  Also in the large energy setting the $\nabla^{-1}(A A)$ term also becomes problematic, indeed it can 
cause breakdown of uniqueness in the Coulomb gauge (though this issue does not arise in the Riemann surface case $m=2$). Because of these difficulties, we shall choose
a different gauge, which we shall call the \emph{canonical heat-temporal gauge} or \emph{caloric gauge}, which will be described later in this section and is constructed using the heat flow for the harmonic map equation; it is the covariant formulation of the
microlocal gauge used in \cite{tao:wavemap}, \cite{tao:wavemap2}, \cite{kr:wavemap}, \cite{tataru:wave3}, and
places $A$ in the schematic form $A = \pi(\psi, \nabla_x^{-1} \psi)$, where $\pi$ is a (nonlinear) paraproduct, defined using the heat flow for the harmonic map equation, which effectively constrains the
$\nabla_x^{-1} \psi$ factor to have higher frequency than the $\psi$ factor.  This is similar in strength to
the Coulomb gauge $\nabla_x^{-1}(\psi \psi)$, but mostly eliminates the ``high-high interactions'' where two high frequency
components of $\psi$ interact to create a low frequency component of $A$; this is precisely the term for
which the $\nabla_x^{-1}$ operator causes difficulties.  Also we do not encounter the $\nabla^{-1}(A A)$ term at all.

As observed in \cite{nahmod}, \cite{shsw:wavemap}, \cite{krieger:3d}, \cite{krieger:2d},
the derivative field $\psi$ evolves by a nonlinear wave equation.  Indeed, from \eqref{wmeq}
we have
$$ \partial_\alpha \psi_\beta - \partial_\beta \psi_\alpha = -A_\alpha \psi_\beta + A_\beta \psi_\alpha$$
and
$$ \partial^\alpha \psi_\alpha = -A^\alpha \psi_\alpha$$
and hence
\begin{equation}\label{psi-wave}
\begin{split}
\partial^\alpha \partial_\alpha \psi_\beta &=
\partial^\alpha (-A_\alpha \psi_\beta + A_\beta \psi_\alpha) + \partial^\alpha \partial_\beta \psi_\alpha \\
&= \partial^\alpha (A_\beta \psi_\alpha - A_\alpha \psi_\beta) - \partial_\beta (A^\alpha \psi_\alpha).
\end{split}
\end{equation}
Given that $A$ should behave like a quadratic expression of $\psi$, containing an inverse derivative (cf. \eqref{hodge}), we thus see that \eqref{psi-wave} is something like a cubic wave equation for $\psi$.

As is well known (see e.g. \cite{shatah-struwe}), the wave map equation enjoys a \emph{stress energy tensor}
$$
\T_{\alpha \beta} := \langle \partial_\alpha \phi, \partial_\beta \phi \rangle_{h}
- \frac{1}{2} g_{\alpha \beta} \langle \partial^\gamma \phi, \partial_\gamma \phi \rangle_{h},
$$
or in terms of the differentiated fields
\begin{equation}\label{stress-def}
\T_{\alpha \beta} := \langle \psi_\alpha, \psi_\beta \rangle_{\R^m}
- \frac{1}{2} g_{\alpha \beta} \langle \psi^\gamma, \psi_\gamma \rangle_{\R^m},
\end{equation}
where $g$ is the Minkowski metric; observe that this quantity is indepedent of the choice of gauge $e$. 
In particular we can define (in standard co-ordinates $x^0 = t, x^1, x^2$)
the \emph{energy density} $\T_{00}$ by
\begin{equation}\label{t00-def}
 \T_{00}(t,x) :=  \frac{1}{2} |\psi_{t,x}(t,x)|^2.
\end{equation}
Because the domain $\R^{1+2}$ is invariant under translations in spacetime, Noether's
theorem yields the conservation law
\begin{equation}\label{conserv}
 \partial^\alpha \T_{\alpha \beta} = 0
\end{equation}
which can of course also be verified directly from \eqref{wmeq}.  The conservation of the energy 
\begin{equation}\label{energy-alt}
E(\phi[t]) = \frac{1}{2} \int_{\R^2} |\psi_{t,x}|^2\ dx = \int_{\R^2} \T_{00}(t,x)\ dx
\end{equation}
now follows directly from \eqref{conserv}.

\begin{remark}  In the case of a Euclidean target $N = \R^m$ (so $\kappa = 0$) and the standard orthonormal gauge
$e = \partial_{x_1}, \ldots, \partial_{x_m}$, we have $\psi_\alpha = \partial_\alpha \phi$, $A_\alpha = 0$, and the 
wave map equation \eqref{cov} then becomes $\partial^\alpha \partial_\alpha \phi = 0$, while \eqref{psi-wave} becomes $\partial^\alpha \partial_\alpha \psi_\beta = 0$.  Thus in this case the wave maps equation collapses to the free (linear) wave equation.
\end{remark}

\subsection{The heat flow equation}

We now leave wave maps for the moment and turn to the heat flow equation for harmonic maps (or \emph{heat flow} for short).  This equation can be written covariantly as
\begin{equation}\label{cov-heat}
 \partial_s \phi = (\phi^*\nabla)_k \partial_k \phi
\end{equation}
where $\phi(s,x) \in N$ is a smooth map defined for some interval $[0,S) \times \R^2$ for some $0 < S \leq +\infty$;
formally, this is gradient flow for the energy functional $\int_{\R^2} \frac{1}{2} |\nabla_x \phi|_{h}^2\ dx$.
Actually, for our applications we will not consider just a single map $\phi(s,x)$ evolving by the heat flow, but an
entire one-parameter family $\phi(s,t,x)$ of such maps, where $t$ ranges in a compact time interval $I$.
While the $s$ variable plays the role of time in this equation, it is quite distinct from the time variable $t$ appearing
in the wave map equation, and so we shall refer to $s$ as the \emph{heat-temporal} variable and $t$ as the
\emph{wave-temporal} variable to minimize confusion.  We use Latin indices $j,k$ to parameterize the spatial variable $x$
(and write $x$ to denote both indices, thus for instance $A_x$ denotes the vector $(A_1, A_2)$),
and Greek indices $\alpha,\beta$ to parameterize the spacetime variables $t,x$, but shall always treat the heat-temporal
variable $s$ separately (using neither Greek nor Latin indices to parameterize these).  In our discussion of the heat
flow the wave-temporal parameter $t$ plays almost no role, although we will need to study the evolution of time derivatives such as $\partial_t \phi$ in the $s$ direction.

In local co-ordinates $\phi^a$, the equation \eqref{cov-heat} becomes a non-linear heat equation
\begin{equation}\label{wm-coord-heat}
\partial_s \phi^a = \Delta \phi^a + \Gamma^a_{bc}(\phi) \partial_\alpha \phi^b \partial^\alpha \phi^c
\end{equation}
as such it is clear (e.g. from energy methods) that this equation is locally well posed in the Schwartz class, and the flow fixes the value $\phi(\infty)$ of the map at infinity.  Indeed, we have the global well-posedness result

\begin{theorem}[Global existence of heat flow]\label{heatflow-conv}  If $\phi: I \times \R^2 \to N$ is Schwartz and equals $\phi(\infty) \in N$ 
at infinity, then there is 
a unique Schwartz solution $\phi: \R^+ \times I \times \R^2 \to N$ to \eqref{cov-heat} with initial data $\phi(0,t,x) = \phi(t,x)$.  Furthermore, $\phi(s)$ converges to $\phi(\infty)$ in the $C^m(I \times \R^2;N)$ topology
as $s \to +\infty$ for any $m \geq 0$.
\end{theorem}

This result is essentially a special case of the results of 
Eells-Sampson \cite{eells} (see also \cite{sy}, \cite{lt}, \cite{lt2}, \cite{struwe:heat}), and follows primarily from
the negative curvature of $N$.  We will not prove this standard result here, but we mention that 
a key tool in this convergence result is the \emph{Bocher identity}
$$ (\partial_s - \Delta)|\psi_\alpha|^2 = \kappa |\psi_\alpha \wedge \psi_x|^2 - 2 |\D_x \psi_\alpha|^2,$$
which in conjunction with the \emph{diamagnetic inequality} $|\nabla_x |\psi_\alpha|| \leq |\D_x \psi_\alpha|$
and the negative curvature hypothesis $\kappa \leq 0$
implies that the evolution of $|\psi_\alpha|$ is dominated by the heat flow:
\begin{equation}\label{subheat}
 (\partial_s - \Delta) |\psi_\alpha| \leq 0.
\end{equation}

We can also write the heat flow equation in differentiated form by lifting $\phi$ to $(\phi,e): \R^+ \times I \times \R^2 \to \Frame(N)$ as before.  We will however limit the gauge freedom here by imposing the \emph{heat-temporal} gauge 
condition 
\begin{equation}\label{heat}
(\phi^* \nabla)_s e = 0;
\end{equation}
in other words, we require that $e$ propagates in the heat-temporal direction $\partial_s$ by parallel transport.  It is clear geometrically (or from the Picard existence theorem) that for fixed $\phi$, the frame $e(s,t,x)$ at $(s,t,x)$ is now determined completely by the corresponding frame $e(0,t,x)$ at $s=0$.  We still have the gauge freedom \eqref{gauge-def}, but only if the gauge $U = U(t,x)$ is independent of $s$, i.e. it depends only on $t,x$.  We will remove this gauge freedom shortly, however, by selecting a canonical heat-temporal gauge.

Using this gauge $e$ we can introduce the derivative fields $\psi_\alpha, \psi_s, A_\alpha, A_s$ as before, although the heat-temporal gauge condition forces $A_s = 0$ (or equivalently $\D_s = \partial_s$).  
From \eqref{cov-cov} and \eqref{psi-def} we see that \eqref{cov-heat} becomes
$\psi_s = \D_k \psi_k$.  In analogy to \eqref{wmeq}, we can now write the heat flow as the system
\begin{equation}\label{heateq}
\begin{split}
\partial_\alpha \phi &= e \psi_\alpha\\
\partial_s \phi &= e \psi_s \\
(\phi^*\nabla)_\alpha e &= e A_\alpha\\
(\phi^*\nabla)_s e = A_s &= 0 \\
\psi_s &= \D_k \psi_k\\
\D_\alpha \psi_\beta - \D_\beta \psi_\alpha &= 0 \\
\partial_s \psi_\alpha = \D_s \psi_\alpha &= \D_\alpha \psi_s\\
F_{\alpha \beta} &= \kappa \psi_\alpha \wedge \psi_\beta\\
\partial_s A_\alpha = F_{s\alpha} &= \kappa \psi_s \wedge \psi_\alpha.
\end{split}
\end{equation}
We refer to the quantity $\psi_s = \D_k \psi_k$ as the \emph{heat-tension field}; it is the direction that the heat flow evolution in the $s$ variable deforms $\phi$ in (with respect to the orthonormal frame $e$), and should be thought of as
a non-linear version of the Laplacian applied to the map $\phi$.   

Similarly to \eqref{psi-wave}, the derivative fields $\psi_\alpha$ obey a (covariant) heat equation
which is linear (but with potential and magnetic components).
Indeed from \eqref{heateq} and \eqref{curv-def} we see that
\begin{equation}\label{psij-eq}
\begin{split}
\partial_s \psi_\alpha &= \D_\alpha \psi_s \\
&= \D_\alpha \D_k \psi_k \\
&= \D_k \D_\alpha \psi_k + F_{\alpha k} \psi_k \\
&= \D_k \D_k \psi_\alpha + \kappa (\psi_\alpha \wedge \psi_k) \psi_k;
\end{split}
\end{equation}
this equation can be viewed as a linearization of the heat flow equation around $\phi$.
Observe that the heat-tension field $\psi_s$ also obeys the same equation:
\begin{equation}\label{psis-eq}
\partial_s \psi_s =  \D_k \D_k \psi_s + \kappa (\psi_s \wedge \psi_k) \psi_k.
\end{equation}
As a first approximation, one can thus think of $\psi_\alpha(s,t)$ as being a non-linear analogue of
$e^{s\Delta} \psi_\alpha(0,t)$, where $e^{s\Delta}$ is the propagator for the free (linear) heat equation.
Similarly for the heat-tension field $\psi_s$.  As is well known, heat operators can be used as a substitute 
for Littlewood-Paley operators; thus one can think of $\psi_\alpha(s,t)$ heuristically as a nonlinear restriction 
of $\psi_\alpha(0,t)$ to frequencies $|\xi| \lesssim s^{-1/2}$.
The heat-tension field $\psi_s$ is similar, but the presence of the derivative in the formula $\psi_s = \D_k \psi_k$ suggests that this quantity is more localized to the annulus $|\xi| \sim s^{-1/2}$ than the ball $|\xi| \lesssim s^{-1/2}$.

\begin{remark}
The assumption $\kappa < 0$ will ensure that the equations \eqref{psij-eq}, \eqref{psis-eq} will not
blow up as $s \to +\infty$; this is basically because of \eqref{subheat}.
\end{remark}

As mentioned above, we still have a gauge freedom to rotate $e$ by an arbitrary gauge $U(t,x)$ independent of $s$.
However, we can remove this freedom by fixing a boundary condition at $s = +\infty$.

\begin{theorem}[Existence of canonical heat-temporal gauge]\label{canonical-heat}
Suppose $\phi: I \times \R^2 \to N$ is Schwartz and equals $\phi(\infty) \in N$ 
at infinity, and let $\phi: \R^+ \times I \times \R^2 \to N$ be the heat flow extension given by
Theorem \ref{heatflow-conv}.  Let $e(\infty)$ be an orthonormal frame in $T_{\phi(\infty)} N$.  Then there is a unique
lift $(\phi,e): \R^+ \times I \times \R^2 \to \Frame(N)$ which is smooth, obeys the heat-temporal 
condition \eqref{heat}, and which converges uniformly to $(\phi(\infty),e(\infty))$ as either $s \to +\infty$ or 
$x \to \infty$.  Furthermore, $A_\alpha(s)$ and $\psi_\alpha(s)$ converge uniformly to zero as $s \to +\infty$.
\end{theorem}

We will not prove this theorem in detail here, but follows from decay estimates on $\psi$ as $s \to +\infty$
which ultimately stem from \eqref{subheat}, which ensure that one can specify boundary data at $s=+\infty$ in a well-posed
manner.

We refer to the gauge $e$ constructed in this manner as the \emph{canonical heat-temporal gauge} or \emph{caloric gauge} associated to $\phi$ (and to $e(\infty)$).  Geometrically, it is constructed by dragging the constant frame $(\phi(\infty),e(\infty))$ back via
parallel transport from $s=+\infty$ to $s=0$ by reversing the heat flow; it can also be viewed as the unique
solution to \eqref{heateq} with boundary conditions
\begin{equation}\label{canonical-gauge}
\begin{split}
(\phi,e) &= (\phi(\infty),e(\infty)) \hbox{ when } s = +\infty\\
(\psi_\alpha,A_\alpha) &= 0 \hbox{ when } s = +\infty.
\end{split}
\end{equation}
Note that we are specifying data at both $s=0$ and $s=+\infty$.  As we shall explain shortly, this gauge is 
the analogue of the microlocal gauge in \cite{tao:wavemap}, \cite{tao:wavemap2}, \cite{kr:wavemap}, \cite{tataru:wave3} 
in the differentiated setting, and will be used here in place of the Coulomb gauge that was used in
\cite{nahmod}, \cite{shsw:wavemap}, \cite{krieger:3d}, \cite{krieger:2d}, for reasons which will also be discussed shortly.

Note that we can now compute $A$ from $\psi$ via the fundamental theorem of calculus and \eqref{heateq}:
\begin{equation}\label{A-fundamental}
\begin{split}
A_\alpha(s,t,x) &= - \int_s^{+\infty} \partial_s A_\alpha(s',t)\ ds' \\
&=  -\kappa \int_s^{+\infty} \psi_s(s',t) \wedge \psi_\alpha(s',t)\ ds'\\
&= - \kappa \int_s^{+\infty} (\D_k \psi_k(s',t)) \wedge \psi_\alpha(s',t)\ ds',
\end{split}
\end{equation}
and in particular we can compute the connection $A_\alpha(t,x)$ on the boundary $s=0$ by the formula
\begin{equation}\label{A-origin}
A_\alpha(0,t,x) = - \kappa \int_0^{+\infty} (\D_k \psi_k(s,t,x)) \wedge \psi_\alpha(s,t,x)\ ds.
\end{equation}
If we continue the heuristic that $\psi(s,t,x)$ is a non-linear version of $e^{s\Delta} \psi(0,t,x)$, then the
above expression can then be viewed as a non-linear analogue of a paraproduct
$\pi( \psi_x(0), \nabla_x^{-1} \psi_x(0) )$,
where the paraproduct $\pi$ is a ``low-high'' paraproduct that restricts the frequency of the $\nabla_x^{-1}\psi_x(0)$
factor to be larger or comparable to that of the $\psi_x(0)$ factor.  To see this informally, we restrict $s$
to the dyadic region $2^{-2\k} \leq s \leq 2^{-2\k+1}$, and then heuristically view $\psi_x(s,t)$ as $P_{\leq \k} \psi_x(0,t)$, a Littlewood-Paley type projection to frequencies $\leq 2^\k$, while $\psi_s(s,t) \D_k \psi_k(s,t)$ is roughly like $2^\k P_\k \psi_x(0,t)$,
a projection to frequencies $\sim 2^\k$.  The expression \eqref{A-origin} is then heuristically like $\sum_\k 
(P_{\leq \k} \psi_x(0,t)) (2^{-\k} P_\k \psi_x(0))$, whence the claim.  
This paraproduct $A = \pi( \psi_x(0), \nabla_x^{-1} \psi_x(0) )$ is favorable as
the inverse derivative $\nabla_x^{-1}$ is always on the highest frequency factor; this should be compared with
the Coulomb gauge, which gives a formula for $A$ of the form $A = \nabla_x^{-1}(\psi_x \psi_x)$ (see \eqref{hodge}), 
which is of similar strength to the paraproduct in \eqref{A-origin} except for high-high interactions of $\psi$, for which
the inverse derivative is in an unfavorable location; this is not particularly harmful in high dimensions (as the arguments
in \cite{nahmod}, \cite{shsw:wavemap} show) but begins
to cause significant problems in two and three dimensions.  In \cite{krieger:3d}, \cite{krieger:2d} the difficulties
arising from this unfavorable placement of the inverse derivative are overcome by using the overdetermined nature of the wave map system to extract additional cancellations from the expression $\nabla_x^{-1}(\psi_x \psi_x)$, and by quite delicate
multilinear estimates.  We will face similar difficulties, but they will be milder because of the more favorable location
of the inverse derivative.

\begin{remark}  We now perform some heuristic manipulations that connect this gauge to the microlocal renormalization
used for the sphere $S^{m-1} \subset \R^m$ in \cite{tao:wavemap}, \cite{tao:wavemap2} (and can also be used to connect to the 
similar procedures used in \cite{kr:wavemap}, \cite{tataru:wave3}).  Here we think of the wave map $\phi$ as a column vector with $m$ entries.  We choose an orthonormal frame $e = (e_1, \ldots, e_k)$ taking values not on the tangent space of the sphere, but rather in the ambient Euclidean space $\R^m$.  In the standard orthonormal frame (which we call $e^0$), the analogues of $\psi_\alpha$ and $A_\alpha$ are
$$ \psi_\alpha := \partial_\alpha \phi; \quad A_\alpha := \partial_\alpha \phi \wedge \phi.$$
In a more general frame $e = U e^0$, where $U$ takes values in the rotation group $O(m)$, we have
$$ \psi_\alpha := U^{-1} \partial_\alpha \phi; \quad A_\alpha := U^{-1} (\partial_\alpha \phi \wedge \phi) U - U^{-1} \partial_\alpha U.$$
The heat-temporal gauge condition $A_s = 0$ thus becomes
$$ \partial_s U = (\partial_s \phi(s) \wedge \phi(s)) U(s)$$
where $\partial_s \phi$ can be defined using the heat flow equation
$$ \partial_s \phi = \Delta \phi + \phi |\nabla \phi|^2.$$
While this equation can blow up for large energies, for small energies it converges to a constant.
In the region $s \sim 2^{-2\k}$, we can then heuristically treat $\phi(s)$ as the low frequency
projection $P_{\leq \k} \phi(0)$, and $\partial_s \phi \approx \Delta \phi$ as the medium frequency
projection $s^{-1} P_\k \phi(0)$.  Thus the equation for $U$ is then approximately
$$ s\partial_s U \approx 
[P_\k \phi(0) \wedge P_{\leq \k} \phi(0)] U(s).$$
Comparing this with the scheme used in \cite{tao:wavemap} (see also \cite{tao:wavemap2}, \cite{kr:wavemap}, \cite{tataru:wave3}) we see that the canonical heat-temporal gauge construction is nothing more than a continuous 
version of the discrete microlocal gauge used in those papers, with the Littlewood-Paley operators in
the ambient Euclidean space being replaced by the more intrinsic operators generated by the harmonic map heat flow.
Note that in the heat flow formulation, the matrices $U$ are automatically orthogonal, whereas in the discrete procedure
in \cite{tao:wavemap}, \cite{tao:wavemap2}, \cite{kr:wavemap}, \cite{tataru:wave3} one only obtains approximate orthogonality, and only when the energy is assumed small.  Indeed for large energies it seems unlikely that one can
replicate this procedure, due to topological obstructions (Schwartz maps from $\R^2$ to $S^{m-1}$ are not necessarily
topologically trivial).
\end{remark}

\begin{remark}  The equation \eqref{A-fundamental} can be used to recover $A$ from $\psi_\alpha$.
One can similarly use \eqref{heateq}, \eqref{A-fundamental}
and the canonical heat gauge condition $\psi_\alpha(\infty,t,x) = 0$
to recover $\psi_\alpha$ from the heat-tension field $\psi_s$:
\begin{equation}\label{psi-gradient}
\begin{split}
\psi_\alpha(s_0,t) &= - \int_{s_0}^{+\infty} \partial_s \psi_\alpha(s,t)\ ds \\
&= -\int_{s_0}^{+\infty} \D_\alpha \psi_s(s,t)\ ds \\
&= - \partial_\alpha \int_{s_0}^\infty \psi_s(s,t)\ ds + \Psi_\alpha(s_0,t)
\end{split}
\end{equation}
where $\Psi$ is the cubic correction term
\begin{equation}\label{Psi-def}
\begin{split}
\Psi_\alpha(s_0,t) &:=  \int_{s_0}^{+\infty} A_\alpha(s) \psi_s(s,t)\ ds\\
&= - \partial_\alpha \int_{s_0}^{+\infty} \psi_s(s)\ ds
+ \kappa \int_{s_0 < s < s'} (\psi_s(s') \wedge \psi_\alpha(s')) \psi_s(s)\ ds' ds.
\end{split}
\end{equation}
Thus one can (at least in principle) recover $\psi_\alpha$ from $\psi_s$ by an iteration process,
and then by \eqref{A-fundamental} recover $A_\alpha$ from $\psi_s$ and $\psi_\alpha$.  Thus one can
view \eqref{heateq} as a scalar non-linear heat equation \eqref{psis-eq} for the heat-tension field $\psi_s$
We will thus often seek to use the above identities to write the other components of $\psi$ and $A$
in terms of this ``dynamic variable'' $\psi_s$ in order to make the system \eqref{waveheat-eq} resemble
 a scalar equation (and thus make more amenable to solution by iterative methods); this corresponds to the 
``dynamic separation'' strategy used in \cite{krieger:3d}, \cite{krieger:2d}.  Indeed one can view the heat-tension field $\psi_s$
as a sort of\footnote{A more precise heuristic is that if $s \sim 2^{-2\k}$, then $\psi_s \sim \frac{1}{s} P_\k \phi$.} continuous Littlewood-Paley resolution of $\phi$, which is the scalar variable in the undifferentiated
wave map equation \eqref{cov}.  One consequence in particular from \eqref{psi-gradient} that the functions $\psi_\alpha$ are 
(modulo a cubic lower order term $\Psi_\alpha$) the gradient of a scalar function $-\int_{s_0}^\infty \psi_s(s)\ ds$.
This fact will be crucial for us, as it means that all of the highest order non-linear expressions 
involving the $\psi_\alpha$ contain null form structure (the lower order terms will turn out to be
manageable via Strichartz estimates and will not require null structure).  In particular we will seek to 
recover null forms such as the trilinear expression $\phi_1 \partial^\alpha \phi_2 \partial_\alpha \phi_3$, 
preferably with $\phi_1$ having ``higher frequency'' than $\phi_2$ or $\phi_3$, in the 
spirit of \cite{tao:wavemap2}; in our heat-flow formulation the $s$ parameter will be a proxy for this concept
of frequency\footnote{More precisely, the quantity $s^{-1/2}$ corresponds to the Fourier notion of frequency magnitude,
as can be seen by inspecting the Fourier symbol of the heat operator $e^{s\Delta}$; thus larger values of $s$
correspond to lower frequencies and vice versa.}.  A similar strategy was carried out in \cite{krieger:3d}, \cite{krieger:2d}, starting not from the heat flow but instead from the equation $\D_\alpha \psi_\beta - \D_\beta \psi_\alpha = 0$ to conclude that $\psi$ has small curl and hence (by Hodge theory) is approximately a gradient.
However the use of Hodge theory once again introduces inverse derivatives $\nabla_x^{-1}$ into the equation, which 
makes the non-linearities more difficult to control than would otherwise be necessary, and so we shall rely on \eqref{psi-gradient} instead of Hodge theory to write $\psi$ as an approximate gradient.
\end{remark}

\subsection{Concatenating the wave map and heat flow equations}

From Theorems \ref{heatflow-conv} and \ref{canonical-heat} we see that the system \eqref{heateq} can be solved, and 
placed in the canonical heat-temporal gauge \eqref{canonical-gauge}, 
given any Schwartz map $\phi: I \times \R^2 \to N$ (as well as an orthonormal frame $e(\infty)$ in $T_{\phi(\infty)} N$.  In particular we can solve this equation for any Schwartz solution to the
wave map equation \eqref{wmeq}.  This gives rise to the fields $\phi$, $e$, $A_\alpha$, $\psi_\alpha$, $A_s$, $\psi_s$ 
on $\R^+ \times I \times \R^2$ solving the system of equations and boundary conditions
\begin{equation}\label{waveheat-eq}
\begin{split}
\partial_\alpha \phi &= e\psi_\alpha \\
\partial_s \phi &= e \psi_s\\
(\phi^*\nabla)_\alpha e &= e A_\alpha\\
(\phi^*\nabla)_s e = A_s &= 0 \\
\D^\alpha \psi_\alpha &= 0 \hbox{ when } s = 0\\
\psi_s &= \D_k \psi_k\\
\D_\alpha \psi_\beta - \D_\beta \psi_\alpha &= 0 \\
\partial_s \psi_\alpha = \D_s \psi_\alpha &= \D_\alpha \psi_s\\
F_{\alpha \beta} &= \kappa \psi_\alpha \wedge \psi_\beta\\
\partial_s A_\alpha = F_{s\alpha} &= \kappa \psi_s \wedge \psi_\alpha\\
(\phi,e) &= (\phi(\infty),e(\infty)) \hbox{ when } s = +\infty\\
(\psi_\alpha,A_\alpha) &= 0 \hbox{ when } s = +\infty.
\end{split}
\end{equation}
Of course, many of these equations are redundant, and this system is highly overdetermined, nevertheless it will
be convenient to retain all of these equations as we will need the full structure of this system at various points
in the argument.

\begin{remark}  If the target $N$ were Euclidean space $N = \R^m$ with the standard frame, then we have $\psi_\alpha = \partial_\alpha \phi$, $\psi_s = \Delta \phi$, $A_\alpha = A_s = 0$, and the heat flow equation \eqref{cov-heat} then becomes $\partial_s \phi = \Delta \phi$, while \eqref{psij-eq} becomes 
$\partial_s \psi_\alpha = \Delta \psi_\alpha$.  Thus in this case the heat flow collapses to the free (linear)
heat equation.  It is clear that \eqref{waveheat-eq} similarly collapses to an evolution which is the wave equation
in the $x,t$ variables and the heat equation in the $x,s$ variables, with the two flows commuting with each other.
For non-Euclidean targets, the same statement is true to top order; thus all fields evolve in the $s$ direction by
a nonlinear heat equation and in the $t$ direction by a nonlinear wave equation.
\end{remark}

\begin{remark} One can use the transitive group action of $O(1,m)$ of $\Frame(N)$ to place $\phi(\infty)$ 
and $e(\infty)$ wherever one pleases.  However there will be no need for us to perform such a normalization.  
\end{remark}

\begin{remark}
We caution the reader that the heat-temporal variable $s$ scales like twice the dimension of space, in contrast
to the wave-temporal variable $t$, which scales like one dimension of space.  Indeed, the scale invariance
of \eqref{waveheat-eq} is given by 
\begin{equation}\label{heat-scale}
\begin{split}
(\phi,e)(s,t,x) &\leftarrow (\phi,e)(\frac{s}{\lambda^2}, \frac{t}{\lambda}, \frac{x}{\lambda})\\
(\psi_\alpha,A_\alpha)(s,t,x) &\leftarrow \frac{1}{\lambda} (\psi_\alpha,A_\alpha)(\frac{s}{\lambda^2}, \frac{t}{\lambda}, \frac{x}{\lambda})\\
\psi_s(s,t,x) &\leftarrow \frac{1}{\lambda^2} \psi_s(\frac{s}{\lambda^2}, \frac{t}{\lambda}, \frac{x}{\lambda})\\
F_{\alpha\beta}(s,t,x) &\leftarrow \frac{1}{\lambda^2} F_{\alpha \beta}(\frac{s}{\lambda^2}, \frac{t}{\lambda}, \frac{x}{\lambda})\\
F_{s\alpha}(s,t,x) &\leftarrow \frac{1}{\lambda^3} F_{s\alpha}(\frac{s}{\lambda^2}, \frac{t}{\lambda}, \frac{x}{\lambda}).
\end{split}
\end{equation}
\end{remark}

Note that the wave map equation $\D^\alpha \psi_\alpha = 0$ only holds
at the boundary $s=0$ of the slab $\R^+ \times I \times \R^2$.  Unfortunately, this equation does not commute with the heat flow equation (even in the completely integrable one-dimensional case $n=1$), so we do not obtain this equation in the interior of the slab.  However, we can assert that the wave map equation holds ``approximately'' in the
interior region.  Indeed, if we define the \emph{wave-tension field} $u$ by
$u := \D^\alpha \psi_\alpha$, then we have (from various equations in
\eqref{waveheat-eq} and covariant versions of the Leibnitz rule)
\begin{equation}\label{u-heat}
\begin{split}
\partial_s u &= \D_s \D^\alpha \psi_\alpha \\
&= F_{s\alpha} \psi^\alpha + \D^\alpha \D_s \psi_\alpha \\
&= \kappa(\psi_s \wedge \psi_\alpha) \psi^\alpha + \D^\alpha \D_\alpha \psi_s \\
&= \kappa(\D_k \psi_k \wedge \psi_\alpha) \psi^\alpha + \D^\alpha \D_\alpha \D_k \psi_k \\
&= \kappa(\D_k \psi_k \wedge \psi_\alpha) \psi^\alpha + 
\D^\alpha (F_{\alpha k} \psi_k) + \D^\alpha \D_k \D_\alpha \psi_k \\
&= \kappa[ (\D_k \psi_k \wedge \psi_\alpha) \psi^\alpha + 
\D^\alpha ((\psi_\alpha \wedge \psi_k) \psi_k)] + \D^\alpha \D_k \D_k \psi_\alpha \\
&= \kappa [(\D_k \psi_k \wedge \psi_\alpha) \psi^\alpha + 
(\D^\alpha \psi_\alpha \wedge \psi_k) \psi_k
+ (\psi_\alpha \wedge \D^\alpha \psi_k) \psi_k
+ (\psi_\alpha \wedge \psi_k) \D^\alpha \psi_k
] \\
&\quad + \F_{\alpha k} \D_k \psi^\alpha + \D_k (\F_{\alpha k} \psi^\alpha)
+ \D_k \D_k \D^\alpha \psi_\alpha \\
&=  \kappa [(\D_k \psi_k \wedge \psi_\alpha) \psi^\alpha + 
(u \wedge \psi_k) \psi_k
+ (\psi_\alpha \wedge \D_k \psi^\alpha) \psi_k
+ (\psi_\alpha \wedge \psi_k) \D^\alpha \psi_k\\
&\quad
+ 2(\psi_\alpha \wedge \psi_k) \D_k \psi^\alpha
+ (\D_k \psi_\alpha \wedge \psi_k) \psi^\alpha
+ (\psi_\alpha \wedge \D_k \psi_k) \psi^\alpha
] + \D_k \D_k u \\
&= \D_k \D_k u + \kappa (u \wedge \psi_k) \psi_k
+ \kappa [(\psi_\alpha \wedge \D_k \psi^\alpha) \psi_k
+ 3(\psi_\alpha \wedge \psi_k) \D_k \psi^\alpha
+ (\D_k \psi_\alpha \wedge \psi_k) \psi^\alpha]\\
&= \D_k \D_k u + \kappa (u \wedge \psi_k) \psi_k + 4 \kappa (\psi_\alpha \wedge \psi_k) \D_k \psi^\alpha.
\end{split}
\end{equation}
Thus the wave-tension field $u$ propogates in the $s$ direction by a covariant, inhomogeneous heat equation.  Also, 
from \eqref{waveheat-eq} we have $u = 0$ on the boundary $s = 0$.  Thus we expect $u$ to be somewhat ``small'' for 
$s>0$ as well.

Just as the wave-tension field obeys a heat equation in $s$, the heat-tension field $\psi_s = \D_k \psi_k$ obeys
a wave equation in $t$:
\begin{equation}\label{psis-wave}
\begin{split}
\partial^\alpha \partial_\alpha \psi_s &=
\partial^\alpha (\D_\alpha \psi_s - A_\alpha \psi_s)\\
&= \partial^\alpha \partial_s \psi_\alpha - \partial^\alpha (A_\alpha \psi_s) \\
&= \partial_s (\D^\alpha - A^\alpha) \psi_\alpha - \partial^\alpha (A_\alpha \psi_s) \\
&= \partial_s u - \partial_s (A^\alpha \psi_\alpha) - \partial^\alpha (A_\alpha \psi_s).
\end{split}
\end{equation}

One can use \eqref{u-heat} and \eqref{psis-wave} to design a (local-in-time) iteration scheme for solving
\eqref{waveheat-eq} in the  large energy setting.  The rigourous details of this scheme will appear elsewhere, but the informal ideas are as follows.  We first use \eqref{u-heat} propagate the
wave-tension field $u$ 
in the $s$ direction (using values of $\psi_{t,x}$ and $A_{t,x}$ obtained from previous iterates), using the boundary condition $u=0$ at $s=0$.  Then one uses \eqref{psis-wave}
to propagate the heat-tension field $\psi_s$ forward in the $t$ direction (using the value of $u$ just obtained, 
and values of $A_\alpha$
and $\psi_s$ obtained from previous iterates), using the initial data of $\psi_s$ at the initial time $t=t_0$. 
 Finally, one uses \eqref{psi-gradient} and \eqref{A-fundamental} to compute $\psi_{t,x}$ and $A_{t,x}$ from $\psi_s$ 
(and using values of $\psi_{t,x}$ obtained from previous iterates), using the boundary condition $\psi_{t,x}=A_{t,x} = 0$ at $s=+\infty$.  This scheme may appear convoluted, but it is necessary to do this in the large energy setting in order
that the iterative algorithms to compute each of the fields $\psi_{t,x}$, $A_{t,x}$, $\psi_s$ eventually involve some integration in time (so that we can take advantage of time localization).  In the small energy case there are simpler schemes available\footnote{For instance, one can propogate $\psi_{t,x}$ forward in $t$ on the boundary $s=0$ using \eqref{psi-wave}, then propagate in $s$ using \eqref{psij-eq}, and then
recover $A$ using \eqref{A-fundamental}.  One can pursue these ideas to obtain a new proof of the small energy regularity results, at least in the constant negative curvature case, which is close in spirit to the arguments in \cite{tao:wavemap2} and particularly in \cite{tataru:wave3}, but we will not do so here.  One could also abandon the caloric gauge, and use the Coulomb gauge instead; for instance, the arguments in \cite{krieger:2d} can (in principle) be modified to obtain an iteration scheme when the target is $\H^2$.  However the $\nabla^{-1}$ factors arising from this gauge would require one to prove far more delicate multilinear estimates, and in particular it becomes even more difficult to take advantage of the time localization than in our current argument.  Also the Coulomb gauge has some uniqueness problems in the large energy non-abelian ($m>2$) case, being in some sense even more non-local than the caloric gauge.}, but the convergence of such schemes in the large energy case becomes problematic because of the possibility of an infinite number of iterations backwards and forwards in the $s$ variable (which cannot be made to converge by localizing $t$).

\begin{remark} An overly simplified model of the system \eqref{waveheat-eq} can be obtained by assuming the
approximation $\psi_\alpha \approx \partial_\alpha \Phi$ for some function $\Phi$ (cf. \eqref{psi-gradient}, \eqref{zero-torsion}), and assuming that the nonlinear heat
flow behaves like the linear heat equation (and in particular behaves somewhat like a family of Littlewood-Paley operators); then the equation for $\Phi$ on the boundary $s=0$ is roughly of the form
$$ \Box \Phi \approx \sum_{\k_1, \k_2, \k_3: \k_1 \geq \min(\k_2,\k_3)} \O((P_{\k_1} \Phi) (\partial_\alpha P_{\k_2} \Phi)
(\partial^\alpha P_{\k_3} \Phi));$$
thus $\Phi$ evolves by a nonlinear wave equation with non-linearity of the form $\Phi \partial_\alpha \Phi \partial^\alpha \Phi$, but with at least one of the derivatives falling on a low frequency term.  This is essentially the type of renormalized equation obtained for the wave map flow when the target is a sphere \cite{tao:wavemap2} or an isometrically embedded manifold in Euclidean space \cite{tataru:wave3}, and it is a good heuristic to keep in mind for the local theory
(including the large data perturbation theory).  For comparison, the corresponding heuristic equation for the
Coulomb gauge formulation of the wave map equation would look something like
$$ \Box \Phi \approx \nabla_x^{-1}[\nabla_x \Phi (\partial_\alpha \Phi)] \partial^\alpha \Phi$$
(see \cite{shsw:wavemap}, \cite{nahmod}, \cite{krieger:3d} and especially \cite{krieger:2d}), which is a similar
type of equation but has an additional divergence arising from ``high-high'' interactions inside the $\nabla_x^{-1}$ term that must be controlled\footnote{There is an additional null structure on the above trilinear form of ``$Q_{ij}$ type'' which can be used to handle this divergence, see \cite{krieger:2d}.  However, in the caloric gauge we neither have the divergence, nor the null structure required to deal with that divergence, and so the analytical treatment is somewhat simpler in this gauge.}.
\end{remark}

\section{Asymptotic local self-similarity}\label{alss-sec}

We now discuss the asymptotic self-similarity of wave maps that concentrate at a point.  Our arguments here
are based entirely on an analysis of the stress energy tensor $\T_{\alpha \beta}$, and in particular are 
independent of the choice of gauge or on the curvature properties of the manifold.

Based in prior 
experience with large data critical non-linear wave equation regularity results (see e.g. \cite{shatah-struwe} for 
some examples), one might hope (in the negative curvature case) to establish 
an energy non-concentration result, for instance showing that
\begin{equation}\label{energy-concentrate}
 \lim_{t \to 0^-} \int_{|x| < |t|} \T_{00}(t,x)\ dx = 0
\end{equation}
whenever $\psi$ is a Schwartz derivative wave map on $[-1,0) \times \R^2$; this for instance would be sufficient
to establish global regularity for large energy wave maps by combining such an energy non-concentration result with a small energy regularity result (e.g. the one in \cite{krieger:2d}) and exploiting a finite speed of propagation argument\footnote{Actually, there is a slight difficulty in truncating the wave map properly to take advantage of finite speed of propagation, as one has to take a little care to ensure that the truncated wave map still has small energy and obeys the required compatibility conditions; we will not discuss these (minor) technicalities here however.  See \cite{tataru:wave3}, \cite{tataru:survey} for some closely related issues.}.  

As is well known, one can hope to obtain \eqref{energy-concentrate}, at least for certain components of the energy
density $\T_{00}$, by contracting the stress energy tensor $\T_{\alpha \beta}$ against some well-chosen vector field\footnote{We shall abuse notation and identify vector fields $X = X^\beta$ with their corresponding first order
differential operators $X = X^\beta \partial_\beta$ without further comment.}
$X^\beta$ in spacetime\footnote{In related equations such as the semilinear wave equation, one often adds lower order correction terms to $\T_{\alpha \beta} X^\beta$; however there are no such lower-order terms
available for wave maps.} and then applying Stokes' theorem on a truncated backwards light cone such as
$\{ |x| < |t|; -T < t < -\eps \}$; the point is that \eqref{conserv} allows us to compute the divergence of $\T_{\alpha \beta}X^{\beta}$ as\footnote{One can of course exploit the symmetry of $\T_{\alpha \beta}$ to replace $\partial^\alpha X^\beta$ by the more symmetric \emph{deformation tensor} $\pi^{\alpha \beta} := \frac{1}{2}(\partial^\alpha X^\beta + \partial^\beta X^\alpha)$, but we shall not need to do so here.}
\begin{equation}\label{deformation}
 \partial^\alpha (\T_{\alpha \beta} X^{\beta}) = \T_{\alpha \beta} \partial^\alpha X^\beta.
\end{equation}
The boundary terms on the cone can be handled by a standard flux decay argument and can
be considered negligible as a first approximation.  The strategy is then to choose $X^\beta$ so that the expression
in \eqref{deformation} (as well as the boundary terms $\T_{0\beta} X^\beta$) consists of terms which are either 
non-negative or small.

In the case of spherical symmetry or equivariant symmetry, one can demonstrate energy decay away from the time axis $\{x=0\}$, not directly by the above strategy, but rather by taking advantage of the symmetry assumptions to eliminate
the angular components of the stress-energy tensor, thus reducing \eqref{conserv} to what is essentially a set of transport equations for the remaining components of the stress-energy tensor in the null directions $\partial_t \pm \partial_r$.  In particular it is relatively
easy (via a Gronwall inequality argument, see e.g. \cite{shatah-struwe}) to show in these special cases that
\begin{equation}\label{time-axis-concentration}
  \lim_{t \to 0^-} \int_{\lambda |t| < |x| < |t|} \T_{00}(t,x)\ dx = 0 \hbox{ for any } 0 < \lambda < 1.
\end{equation}
If one then applies \eqref{deformation} to vector fields such as $X = r\partial_r = x \cdot \nabla_x$ and applies
the above estimate, one can easily obtain (time-averaged) decay of the time component of the energy:
\begin{equation}\label{time-decay}
  \lim_{T \to 0^-} \frac{1}{T} \int_{-T}^0 \int_{|x| < |t|} |\psi_0(t,x)|^2\ dx dt = 0.
\end{equation}
From this and the negative curvature of $N$ we can then obtain decay of the spatial component of the energy also.

In the absence of symmetry, one cannot hope to apply this type of argument directly.  First of all, the stress-energy
tensor no longer propagates solely in the null directions $\partial_t \pm \partial_r$, but can now propogate in any timelike or lightlike direction, which seems to frustrate any attempt to use Gronwall's inequality to obtain much decay, even very close to the light cone.  A more serious objection, however, arises from the Lorentz invariance of
the wave map equation once the symmetry assumptions are removed.  Indeed, if one could prove the decay estimate
\eqref{time-axis-concentration} away from the time axis for arbitrary Schwartz wave maps, then by Lorentz invariance
one could also prove a similar decay estimate away from any other timelike ray emenating backwards in time from
the spacetime origin $(0,0)$.  Using two such disjoint timelike rays, it is then an easy matter to obtain
\eqref{energy-concentrate}.  Thus there is no advantage in excluding the time axis in proving \eqref{time-axis-concentration} for general wave maps; this estimate is as difficult as the original 
estimate \eqref{energy-concentrate}.  In particular, one cannot hope to obtain an estimate such as \eqref{time-axis-concentration} without\footnote{For instance, the numerical work in \cite{bct} strongly suggests that when the target is $S^2$ and the initial data is equivariant, then the energy can concentrate along the time axis leading to blowup.
Applying a Lorentz transformation, one can then obtain non-equivariant initial data which concentrates
on another timelike ray, and in particular \eqref{time-axis-concentration} fails for this choice of data.}
relying more heavily on the negative curvature of $N$; this in particular rules out the
possibility of proving \eqref{time-axis-concentration} purely from an analysis of the stress-energy tensor.  Similar
considerations hold for \eqref{time-decay}, basically because the vector field $\partial_t$ is not Lorentz-invariant.

However, a closer inspection of \eqref{time-decay} shows that, in light of \eqref{time-axis-concentration}, the content
of \eqref{time-decay} is only new when $(t,x)$ is close to the time axis; thus we only need the relevant vector field
to point in the direction of $\partial_t$ near the time axis.  This now allows for the possibility of a Lorentz-invariant
estimate, for instance using the scaling vector field $S := x^\alpha \partial_\alpha = t \partial_t + x \cdot \nabla_x$, or more generally $S / \rho^k$, where $\rho := \sqrt{- x^\alpha x_\alpha} = (|t|^2 - |x|^2)^{1/2}$
and $k \in \R$ is some parameter.  For \emph{harmonic} maps in three dimensions $\R^3$ it is well known (see e.g. \cite{schoen}) that the choice
$k=1$ (which would give the radial vector field $\partial_r$ in three dimensions) yields a useful monotonicity
formula of Pohozaev type; in light of the algebraic analogy between Euclidean space $\R^3$ and Minkowski space
$\R^{1+2}$ it is then natural to try the same thing for wave maps on $\R^{1+2}$.  This indeed works, and
gives the following decay estimate:

\begin{proposition}[Asymptotic local self-similarity]\label{alss}  Let $\psi$ be a Schwartz derivative wave map on $[-1,0) \times \R^2$ (not assumed to have any symmetry assumptions).  Then
\begin{equation}\label{scaled-decay}
  \lim_{T \to 0^-} \frac{1}{|\log |T||} \int_{-1}^{-T} \int_{|x| < |t|} |\frac{1}{t} \psi_S(t,x)|^2\ dx \frac{dt}{t} = 0,
\end{equation}
where $\psi_S$ is the component of $\psi$ in the scaling direction $S$, thus
$$ \frac{1}{t} \psi_S := \frac{1}{t} S^\alpha \psi_\alpha = \psi_0 + \frac{x_j}{t} \psi_j.$$
\end{proposition}

This estimate should be compared with \eqref{time-decay}.  On one hand, it replaces the vector field $\partial_t$
by the variant $\frac{1}{t} S = \partial_t + \frac{x}{t} \cdot \nabla_x$, which is
(conformally) invariant under mild Lorentz transformations (i.e. Lorentz transformations with bounded coefficients\footnote{We do not expect our estimates to be invariant under extreme Lorentz transformations since 
such transformations tend to increase the energy substantially.}) and applies for general
wave maps.  On the other hand, it requires somewhat more averaging in time than is present in \eqref{time-decay}.
Note that conservation of energy only allows us to say that the expression in the left-hand side of \eqref{scaled-decay} is bounded, but not that it decays to zero.  

We now give a proof of the above proposition.

\begin{proof}[Proof of Proposition \ref{alss}]
Let $\psi$ be a smooth derivative wave map on $[-1,0) \times \R^2$.  We now elaborate the stress energy analysis
begun in \eqref{deformation}.  For any two times $-1 \leq t_2 < t_1 < 0$, define the truncated solid cone
$$ K_{[t_2,t_1]} := \{ (t,x): t_2 \leq t \leq t_1; |x| \leq |t| \}.$$
This truncated cone has the boundary $\partial K_{t_2} \cup \partial K_{(t_2,t_1)} \cup \partial K_{t_1}$, where
$\partial K_t$ is the disk 
$$ \partial K_t := \{ (t,x): |x| \leq |t| \} \subset \{t\} \times \R^2$$
and $\partial K_{(t_2,t_1)}$ is the truncated light cone
$$ \partial K_{(t_2,t_1)} := \{ (t,x): t_2 < t < t_1; |x| = |t| \}.$$
The surfaces $\partial K_t$ can be given Lebesgue measure $dx$ (which co-incides with the measure induced from the ambient Minkowski metric).  The surfaces $\partial K_{(t_2,t_1)}$ are null and thus have no canonical measure, however we will give them the artificial measure $d\sigma$, defined by
$$ \int_{\partial K_{(t_2,t_1)}} u\ d\sigma := \int_{t_2}^{t_1} \int_{S^1} u(t, |t|\omega)\ |t| d\omega dt$$
for all test functions $u$, where $d\omega$ is arclength measure on the circle $S_1$.  From Stokes' theorem we thus have
$$ -\int_{K_{[t_2,t_1]}} \partial^\alpha P_\alpha(t,x)\ dt dx
= \int_{\partial K_{t_1}} P_0(t_2,x)\ dx + \int_{\partial K_{(t_2,t_1)}} P_L\ d\sigma
- \int_{\partial K_{t_2}} P_0(t_1,x)\ dx,$$
for any one-form $P_\alpha$ smooth on $K_{[t_2,t_1]}$, where
$$ P_L(t,x) := P_0(t,x) - \frac{x_j}{|x|} P_j(t,x)$$
is the $L := \partial_t - \partial_r$ component of $P$.  In particular, from \eqref{deformation} we have
\begin{equation}\label{general-stokes}
-\int_{K_{[t_2,t_1]}} \T_{\alpha \beta} \partial^\alpha X^\beta\ dx dt
= \int_{\partial K_{t_1}} \T_{0\beta} X^\beta\ dx + \int_{\partial K_{(t_2,t_1)}} 
\T_{L \beta} X^\beta\ d\sigma
- \int_{\partial K_{t_2}} \T_{0\beta} X^\beta\ dx
\end{equation}
for any vector field $X^\beta$ smooth on $K_{[t_2,t_1]}$.

We can now derive the standard facts of energy monotonicity and flux decay by choosing the vector field
$X = \partial_t$ in \eqref{general-stokes}, obtaining the energy identity
\begin{equation}\label{energy-identity}
 0 = \int_{\partial K_{t_1}} \T_{00}\ dx + \int_{\partial K_{(t_2,t_1)}} 
\T_{L0}\ d\sigma
- \int_{\partial K_{t_2}} \T_{00}\ dx.
\end{equation}
From \eqref{stress-def} one easily verifies that
$$ T_{L0} = \frac{1}{2} |\psi_L|_{N}^2 + \frac{1}{2} |\psi_\nabb|_{N}^2 \geq 0$$
on $\partial K_{(t_2,t_1)}$, where $\nabb := \frac{x_1}{|x|} \partial_{x_2} - \frac{x_2}{|x|} \partial_{x_1}$
is the spatial angular derivative, thus $\psi_\nabb = \frac{x_1}{|x|} \psi_2 - \frac{x_2}{|x|} \psi_1$.  
In particular, we have \emph{energy monotonicity}
$$ \int_{\partial K_{t_1}} \T_{00}\ dx \leq \int_{\partial K_{t_2}} \T_{00}\ dx \hbox{ for all } -1 \leq t_2 \leq t_1 < 0;$$
if we then set $E_0 := \int_{\partial K_{-1}} \T_{00}\ dx$, then we have
\begin{equation}\label{energy-bound}
 \int_{\partial K_{t}} \T_{00}\ dx \leq E_0 \hbox{ for all } -1 \leq t < 0.
\end{equation}
From this, \eqref{energy-identity}, and monotone convergence we obtain the bounded flux property
$$ \int_{\partial K_{(-1,0)}} \T_{L0}\ d\sigma \leq E_0 < \infty.$$
In particular we have \emph{flux decay}
\begin{equation}\label{flux-decay}
 \lim_{t \to 0^-} \int_{\partial K_{(t,0)}} \T_{L0}\ d\sigma = 0.
\end{equation}

Now let $\lambda > 4$ be a large number, and let $-1/\lambda < t_1 < 0$.  We set $t_2 := \lambda t_1$, and 
let $\eta(t,x) = \eta(t)$ be a time cutoff which 
equals 1 when $t_2/2 < t < 2t_1$, vanishes when $t < t_2$ or $t > t_1$, and is smooth otherwise (so in 
particular $\eta'(t) = O(1/t_2)$ when $t_2 \leq t \leq t_2/2$ and $\eta'(t) = O(1/t_1)$ when $2t_1 \leq t \leq t_1$).
We would now like to apply \eqref{general-stokes} to the vector field
$$ X^\beta := \frac{\eta x^\beta}{\rho} + (\partial^\beta \eta) \rho,$$
where $\rho := \sqrt{-x^\alpha x_\alpha} = \sqrt{|t|^2 - |x|^2}$, as the left-hand side of \eqref{general-stokes} 
will then give something very similar to the expression in \eqref{scaled-decay}.  Unfortunately this vector field 
is singular on the light cone and so \eqref{general-stokes} does not directly apply.  To resolve this problem we shall
mollify the above vector field slightly.  Let $0 < \eps \leq 1$ be a small parameter, and define the modified co-ordinates
$\tilde x^{\alpha}$ and $\tilde \rho$ by
$$\tilde x^0 := x^0 + \eps t_1; \quad \tilde x^j := x^j; \quad \tilde \rho = \sqrt{- \tilde x^{\alpha} \tilde x_\alpha}$$
for $j=1,2$; thus the $\tilde x^\alpha$ co-ordinates are just the $x^\alpha$ co-ordinates but with the origin shifted
upwards to $(-\eps t_1, 0)$.  In particular $\tilde \rho$ is now smooth on all of $K_{[t_2,t_1]}$.  We now apply
\eqref{general-stokes} with the vector field
$$ \tilde X^\beta := \frac{\eta \tilde x^\beta}{\tilde \rho} + (\partial^\beta \eta) \tilde \rho.$$
Using the identities 
$$\tilde \rho^2 = -\tilde x^\gamma x_\gamma; \quad \partial^\alpha \tilde \rho = - \frac{\tilde x^\alpha}{\tilde \rho};  \quad \partial^\alpha \tilde x^\beta = g^{\alpha \beta},$$
where $g$ is the Minkowski metric, we see that
$$ \partial^\alpha \tilde X^\beta = \frac{\eta}{\tilde \rho^3} ( 
\tilde x^\alpha \tilde x^\beta - g^{\alpha \beta} \tilde x^\gamma \tilde x_\gamma ) + \rho \partial^\alpha \partial^\beta \eta.$$
Contracting this against $\T_{\alpha \beta}$ using \eqref{stress-def} (and the fact that $g^{\alpha \beta} g_{\alpha \beta} = \dim(\R^{1+2}) = 3$), we obtain
$$ T_{\alpha \beta} \partial^\alpha \tilde X^\beta = \frac{\eta}{\tilde \rho^3} 
|\tilde x^\alpha \psi_\alpha|^2
+ \rho \T_{\alpha \beta} \partial^\alpha \partial^\beta \eta.$$
Inserting this into \eqref{general-stokes}, and observing that the cutoff $\eta$ eliminates the boundary terms at
$\partial K_{t_1}$, $\partial K_{t_2}$, we obtain
\begin{equation}\label{interesting}
-\int_{K_{[t_2,t_1]}} \frac{\eta}{\tilde \rho^3} 
|\tilde x^\alpha \psi_\alpha|^2 + \rho \T_{\alpha \beta} \partial^\alpha \partial^\beta \eta\ dx dt
= \int_{\partial K_{(t_2,t_1)}}
\T_{L \beta} (\frac{\eta \tilde x^\beta}{\tilde \rho} + (\partial^\beta \eta) \tilde \rho)\ d\sigma.
\end{equation}
Consider the second term on the left-hand side, which simplifies to $\rho \T_{00} \eta''$.  This is supported on the regions $t_2 < t < t_2/2$ and $2t_1 < t < t_1$, and we have the crude bound $\rho \eta'' = O(1/t)$ on those regions.
Thus by the energy bound \eqref{energy-bound}, the contribution of this second term is $O(E_0)$.  Now consider the right-hand side of \eqref{interesting}.  
From \eqref{stress-def} it is easy to show\footnote{This is
best seen using the null frame $L$, $\nabb$, and $\underline L := \partial_t + \partial_r$; indeed we have
$\T_{LL} = |\psi_L|^2$, $\T_{L\nabb} = \langle \psi_L, \psi_\nabb \rangle$, and $\T_{L\underline{L}} =
- \frac{1}{2} |\psi_\nabb|^2$.}
 that $\T_{L\beta} = O(\T_{L0})$ for all
$\beta$.
Also, we have the pointwise bound
$$ |\frac{\eta \tilde x^\beta}{\tilde \rho} + (\partial^\beta \eta) \tilde \rho| \leq C(\lambda,\eps) < \infty$$
on $\partial K_{(t_2,t_1)}$, for some quantity $C(\lambda,\eps)$ depending on $\lambda$ and $\eps$ (the value of $C$ may vary from line to line).  Thus we have
$$ \int_{K_{[t_2,t_1]}} \frac{\eta}{\tilde \rho^3} |\tilde x^\alpha \psi_\alpha|^2\ dx dt
\leq C E_0 + C(\lambda,\eps) \int_{\partial K_{(t_2,t_1)}} \T_{L0}\ d\sigma.$$
Applying flux decay \eqref{flux-decay} we thus have
$$ \limsup_{t_1 \to 0^-}
\int_{K_{[\lambda t_1,t_1]}} \frac{\eta}{\tilde \rho^3} |\tilde x^\alpha \psi_\alpha|^2\ dx dt
\leq C E_0$$
uniformly for all choices of $\lambda$ and $\eps$.  We shall not take advantage of the ability to set $\eps$
to be small, and just take $\eps := 1$.  Then by the triangle inequality we have
$$ |\psi_S| = |x^\alpha \psi_\alpha| \leq |\tilde x^\alpha \psi_\alpha|
+ |t_1| |\psi_0|.$$
On the other hand, using the crude bound $\frac{1}{\tilde \rho^3} \leq \frac{C}{|t|^3}$ and \eqref{energy-bound} 
we see that
$$
\int_{K_{[\lambda t_1,t_1]}} \frac{\eta}{\tilde \rho^3} |t_1|^2 |\psi_0|^2\ dx dt
\leq C E_0$$
uniformly in $\lambda$ and $t_1$, and thus
$$ \limsup_{t_1 \to 0^-} \int_{K_{[\lambda t_1/2, 2t_1]}} \frac{1}{|t|^3}
|\psi_S|^2\ dx dt \leq C E_0.$$
In particular, for any $\lambda$ we see that for $t_1$ sufficiently close to 0 (depending on $\lambda$) we have
$$\int_{K_{[\lambda t_1/2, 2t_1]} }
|\frac{1}{t} \psi_S|^2\ dx \frac{dt}{t} \leq C E_0.$$
Telescoping this, we obtain that
$$\limsup_{T \to 0^-} \frac{1}{|\log |T||} \int_{K_{[-1, T]} }
|\frac{1}{t} \psi_S|^2\ dx \frac{dt}{t} \leq C E_0 / \log \lambda$$
(note that one can replace the interval $[-1,T]$ by $[t_0,T]$ for any fixed $t_0$ without affecting the limit
superior).  Letting $\lambda \to \infty$ we obtain \eqref{scaled-decay}.
\end{proof}

\begin{remark} The estimate works for all targets $N$, as it did not use the hypothesis that $N$ had constant
negative curvature.  It is also completely covariant and did not rely the fixing of a gauge.
\end{remark}

\begin{remark}  Observe that in the above proof we have in fact obtained stronger estimates than \eqref{scaled-decay};
for instance, we can improve the decay of $\frac{1}{t} \psi_S$ near the light cone,
and obtain good estimates on relatively narrow slabs such as $K_{[\lambda t, t]}$ instead of the tall slabs $K_{[-1,T]}$
in \eqref{scaled-decay}.  For instance, it is not difficult to use the above arguments to obtain a sequence of
times $t_k \to 0^-$, a sequence of scales $\lambda_k \to \infty$, and $\eps_k \to 0^+$, such that
$$ 
\int_{K_{[\lambda_k t_k,t_k]}} \frac{|\psi_S|_{N}^2}{(\rho + \eps_k |t|)^3}\ dx dt
\leq \eps_k.$$
Morally speaking, this suggests that the rescaled original wave maps $\phi(t/t_k,x/t_k)$
are becoming locally self-similar as $k \to \infty$.  Unfortunately this is not by itself enough to deduce that these
rescaled wave maps converge (weakly) to a \emph{non-trivial} self-similar wave map, even after other rescalings,
because the above estimates do not prevent the energy from dispersing into multiple points of concentration, or disappearing into the light cone; the basic problem is that control of the $\psi_S$ component of the energy density does not seem to directly control the other components without further structural control on $\phi$, even near the light cone (the presence of angular components in the
stress energy seems to prevent the Gronwall inequality approach in \cite{christ.spherical.wave}, \cite{shatah-struwe}, etc.
from being effective, and naive attempts to exploit the negative curvature of the target $N$ seem to require
more boundedness control on $\phi$ than is currently available.)
\end{remark}

\end{document}